\documentstyle[11pt,amstex,amssymb]{amsart}\textheight=19.8truecm\textwidth=14.1truecm
\begin{document}

\ \vskip 30pt \centerline{\Large\bf A matrix subadditivity
inequality for $f(A+B)$ and $f(A)+f(B)$}

\vskip 15pt \centerline{{\sl \ Jean-Christophe Bourin \,and\,
Mitsuru Uchiyama}} \vskip 5pt \centerline{{\small December 4, 2006}}

\vskip 30pt \noindent {\small {\bf Abstract.} \vskip 5pt In 1999
Ando and Zhan proved a subadditivity inequality for {\it operator}
concave functions. We extend it to {\it all} concave functions:
Given positive semidefinite matrices $A$, $B$ and a non-negative
concave function $f$ on $[0,\infty)$,
$$
\Vert f(A+B)\Vert \le \Vert f(A)+f(B)\Vert
$$

\noindent for all symmetric norms (in particular for all Schatten
$p$-norms). The case $f(t)=\sqrt{t}$ is connected to some
block-matrix inequalities, for instance the operator norm inequality
$$
\left\|
\begin{pmatrix}
A & X^* \\
X & B
\end{pmatrix}
\right\|_{\infty} \le \max\{\, \left\|  \, |A|+|X|\,
\right\|_{\infty} ; \left\|  \, |B|+|X^*|\, \right\|_{\infty}\,\}
$$
for any partitioned Hermitian matrix.


Mathematical subjects classification:   15A60, 47A30}

\vskip 25pt

\vskip 20pt\noindent {\large\bf 1. A subadditivity inequality}

 \vskip 10pt Capital letters $A$, $B,\dots,Z$ mean $n$-by-$n$ complex matrices, or operators on an
  $n$-dimensional Hilbert space ${\cal H}$. If $A$ is positive semidefinite, resp.\ positive definite, we write $A\ge 0$, resp.\ $A>0$.
Recall that a symmetric (or unitarily invariant) norm
$\Vert\cdot\Vert$ satisfies  $\Vert A\Vert = \Vert UAV\Vert$ for all
$A$ and all unitaries $U,\, V$. We will prove:

\vskip 10pt\noindent
 {\bf Theorem 1.1.}  {\it Let $A,\, B\ge 0$ and let  $f:[0,\infty)\longrightarrow[0,\infty)$ be a concave
function. Then, for all symmetric norms,
$$
\Vert f(A+B)\Vert \le \Vert f(A)+f(B)\Vert.
$$
}

\vskip 10pt\noindent For the trace norm Theorem 1.1 is a classical
inequality. In case of the operator norm, Kosem [7] recently gave
a three-line proof! But the general case is much more difficult.
When $f$ is operator concave, Theorem 1.1 has been proved by Ando
and Zhan [1]. Their proof is not elementary and makes use of
integral representations of operator concave functions. By a quite
ingenious process, Kosem [7] derived from Ando-Zhan's result a
related superadditive inequality:

  \vskip 10pt\noindent
   {\bf Theorem 1.2.}  {\it Let $A,\, B\ge 0$ and let
$g:[0,\infty)\longrightarrow[0,\infty)$ be  a convex function with
$g(0)=0$. Then, for all symmetric norms,
$$
\Vert g(A) + g(B)\Vert \le \Vert g(A+B)\Vert.
$$
}

\noindent The special case $g(t)=t^m,\, m=1,2,\dots$ is due to
Bhatia-Kittaneh [4]. The general case has been conjectured by
Aujla and Silva [3].

In this note we first give a simple proof of these two theorems. Our
method is elementary: we only use a simple inequality for operator
convex functions and some basic facts about symmetric norms and
majorization. For background we refer to [9] and references herein.

Next, we consider some inequalities for block-matrices inspired by
the observation that Theorem 1.1 for $f(t)=\sqrt{t}$ can be written
as $\Vert \sqrt{A^2 + B^2}\Vert \le \Vert A + B\Vert$, or
equivalently
$$
\left\|
\begin{pmatrix}
A & 0 \\
B & 0
\end{pmatrix}
\right\| \le \left\| A+B \right\|.
$$
We naturally asked if a similar result holds when the zeros are
replaced by arbitrary positive matrices. We got proofs from T.\
Ando, E.\ Ricard and X.\ Zhan. We thank them for their
collaboration.

\vskip 20pt\noindent {\large\bf 2. \ Proof of Theorems 1.1-1.2 and
related results}

\vskip 10pt First we sketch the simple proof from [8] for Theorem
1.1 in the operator concave case (Ando-Zhan's inequality). Let us
recall some basic facts about operator convex/concave functions on
an interval $[a,b]$. If $g$ is operator convex and $A$ is Hermitian,
$a\ge A \ge b$, then for all subspaces ${\cal S}\subset{\cal H}$,
Davis' Inequality holds for compressions onto ${\cal S}$,
\begin{equation}
g(A_{\cal S}) \le g(A)_{\cal S}.
\end{equation}
Assuming $0\in[a,b]$, $g(0) \le 0$, one can derive Hansen's
Inequality: {\it $Z$ being any contraction,}
$$
g(Z^*AZ) \le Z^*g(A)Z.
$$
Of course, for an operator concave function $f$ on $[a,b]$ with
$f(0) \ge 0$, the reverse inequality holds. For such an $f$ on the
positive half-line and $A,\, B>0$ we then have
$$
f(A) \ge A^{1/2}(A+B)^{-1/2}f(A+B)(A+B)^{-1/2}A^{1/2}
$$
since $Z=(A+B)^{-1/2}A^{1/2}$ is a contraction and $A=Z^*(A+B)Z$.
Similarly
$$
f(B) \ge B^{1/2}(A+B)^{-1/2}f(A+B)(A+B)^{-1/2}B^{1/2}.
$$
Consequently
\begin{equation}
f(A)+f(B)  \ge A^{1/2}\frac{f(A+B)}{(A+B)}A^{1/2} +
B^{1/2}\frac{f(A+B)}{(A+B)}B^{1/2}.
\end{equation}

 Next, the main observation of [8] can be stated as:

 \vskip 10pt\noindent
 {\bf Proposition 2.1.}  {\it Let $A,\, B\ge 0$ and let
$g:[0,\infty)\longrightarrow[0,\infty)$. If $g(t)$ decreases and
$tg(t)$ increases, then for all symmetric norms,
$$
\Vert (A+B)g(A+B)\Vert \le \Vert A^{1/2}g(A+B)A^{1/2} +
B^{1/2}g(A+B)B^{1/2}\Vert.
$$
}

 Combining (2) and Proposition 2.1 with $g(t)=f(t)/t$ yields the
Ando-Zhan Inequality [1]:

 \vskip 10pt\noindent
 {\bf Corollary 2.2.}  {\it Theorem 1.1 holds when $f$ is operator concave. }

\vskip 10pt\noindent This means that the eigenvalues of $f(A+B)$ are
weakly majorised by those of $f(A) + f(B)$. Suppose now that $f$ is
onto, thus $f(0)=0$, $f(\infty)=\infty$ and its inverse function $g$
is convex, increasing. Therefore the eigenvalues of $g(f(A+B))=A+B$
are weakly majorised by those of $g(f(A) + f(B))$. Replacing $A$ and
$B$ by $g(A)$ and $g(B)$ respectively, we get the second Ando-Zhan
Inequality [1]:

\vskip 10pt\noindent
 {\bf Corollary 2.3.}  {\it Let $A,\, B\ge 0$ and let
$g:[0,\infty)\longrightarrow[0,\infty)$ be  a one to one function
whose inverse function is operator concave. Then, for all
symmetric norms,
$$
\Vert g(A) + g(B)\Vert \le \Vert g(A+B)\Vert.
$$
}

\vskip 10pt Now we turn to a quite simple proof of Theorem 1.2. It
suffices to consider Ky Fan $k$-norms $\Vert\cdot\Vert_{k}$. Suppose
that $f$ and $g$ both satisfy Theorem 1.2. Using the triangle
inequality and the fact that $f$ and $g$ are nondecreasing,
\begin{align*}
\Vert (f+g)(A) + (f+g)(B)\Vert_{k}&\le \Vert f(A)+f(B)\Vert_{k} + \Vert g(A)+g(B)\Vert_{k} \\
&\le \Vert f(A+B)\Vert_{k} + \Vert g(A+B)\Vert_{k} = \Vert
(f+g)(A+B)\Vert_{k},
\end{align*}
hence the set of functions satisfying  Theorem 1.2 is a cone. It
is also closed for pointwise convergence. Since any positive
convex function vanishing at $0$ can be approached by a positive
combination of angle functions at $a>0$,
$$
\gamma(t) =\frac{1}{2}\{ |t-a| + t-a\},
$$
it suffices to prove Theorem 1.2 for such a $\gamma$. By Corollary
2.3 it suffices to approach $\gamma$ by functions whose inverses are
operator concave. We take (with $r>0$)
$$
h_r(t)=\frac{1}{2}\{ \sqrt{(t-a)^2 +r} + t -\sqrt{a^2 +r}\},
$$
whose inverse
$$
t-\frac{r/2}{2t+\sqrt{a^2 +r}-a} + \frac{\sqrt{a^2 +r}+a}{2}
$$
is operator concave since $1/t$ is operator convex on the positive
half-line (inequality (1) is then a basic fact of Linear Algebra).
Clearly, as $r\to 0$, $ h_r(t)$ converges uniformly to $\gamma$.

From Theorem 1.2 we can derive Theorem 1.1:

\vskip 10pt\noindent
 {\bf Proof of Theorem 1.1.} It suffices to
prove the theorem for the Ky Fan $k$-norms $\Vert\cdot\Vert_{k}$.
This shows that we may assume  $f(0)=0$. Note that $f$ is
necessarily non-decreasing. Hence, there exists a rank $k$ spectral
projection $E$ for $A+B$, corresponding to the $k$-largest
eigenvalues $\lambda_1(A+B),\dots,\lambda_k(A+B)$ of $A+B$, such
that
$$
\Vert f(A+B)\Vert_{k}=\sum_{j=1}^k \lambda_j(f(A+B))={\rm Tr\,} Ef(A+B)E.
$$
Therefore, using a well-known property of Ky Fan norms, it suffices to show that
$$
{\rm Tr\,} Ef(A+B)E \le {\rm Tr\,} E(f(A)+f(B))E.
$$
This is the same as requiring that
\begin{equation}
{\rm Tr\,} E(g(A)+g(B))E \le {\rm Tr\,} Eg(A+B)E
\end{equation}
for all non-positive convex functions $g$ on $[0,\infty)$ with
$g(0)=0$. Any such function can be approached by a combination of
the type
\begin{equation*}
g(t)=\lambda t + h(t)
\end{equation*}
for a scalar $\lambda <0$ and some non-negative convex function $h$
vanishing at 0. Hence, it suffices to show that (3) holds
 for $h(t)$.
We have
\begin{align*}
{\rm Tr\,} E(h(A)+h(B))E &= \sum_{j=1}^k \lambda_j(E(h(A)+h(B))E) \\
&\le \sum_{j=1}^k \lambda_j(h(A)+h(B)) \\
&\le \sum_{j=1}^k \lambda_j(h(A+B)) \quad {\rm (by\ Theorem\ 1.2)} \\
&= \sum_{j=1}^k \lambda_j(Eh(A+B)E) \\
&= {\rm Tr\,} Eh(A+B)E
\end{align*}
where the second equality follows from the fact that $h$ is
non-decreasing and hence $E$ is also a spectral projection of
$h(A+B)$ corresponding to the $k$ largest eigenvalues.
 \qquad $\Box$

 \vskip 10pt The above proof  is inspired by a part
of the proof of the following result [5]:

 \vskip 10pt \noindent {\bf Theorem 2.4.} {\it Let
$f:[0,\infty)\longrightarrow [0,\infty)$ be a concave function. Let
$A\ge0$  and let $Z$ be expansive. Then, for all symmetric norms,
\begin{equation*}
\Vert f(Z^*AZ) \Vert \le \Vert Z^*f(A)Z\Vert.
\end{equation*}
}

\noindent Here $Z$ expansive means $Z^*Z\ge I$, the identity
operator. Besides such symmetric norms inequalities, there also
exist interesting inequalities involving unitary congruences. For
instance [2]:

\vskip 10pt\noindent {\bf Theorem 2.5.}  {\it Let $A,\, B\ge 0$ and
let  $f:[0,\infty)\longrightarrow[0,\infty)$ be  a concave function.
Then, there exist unitaries $U,\, V$ such that
$$
f(A+B) \le Uf(A)U^* + Vf(B)V^*.
$$
}

\noindent This implies that $ \lambda_{j+k+1}f(A+B) \le
\lambda_{j+1}f(A) + \lambda_{k+1}f(B) $ for all integers $j,\,
k\ge0$.

Combining Theorem 2.5 with Thompson's triangle inequality we get:

\vskip 10pt\noindent {\bf Corollary 2.6.} {\it For any $A$, $B$ and
any non-negative concave function $f$ on $[0,\infty)$,
$$
f(|A+B|) \le Uf(|A|)U^* + Vf(|B|)V^*
$$
for some unitaries $U,$ $V$. }

\vskip 10pt\noindent
 Therefore, we recapture a result from
[8]: the map $X\longrightarrow \Vert f(X) \Vert$ is subbaditive. For
the trace-norm, this is Rotfel'd Theorem.

\vskip 10pt
  A remark may be added
about Theorem 1.1: It can be
 stated for a family $\{A_i\}_{i=1}^m$ of positive operators,
$$
\Vert f(A_1+\dots+A_m) \Vert \le  \Vert f(A_1)+\dots +f(A_m)\Vert.
$$
Indeed, Proposition 2.1 can be stated for a suitable family
$A,\,B,\dots.$

 \vskip 20pt\noindent {\large\bf 3. \ Inequalities for
block-matrices}

\vskip 10pt In Section 1 we noted an inequality involving a
partitioned matrix. The following two theorems are generalizations
due to T.\ Ando, E.\ Ricard and X.\ Zhan (private communications).
The symbol $||\cdot||_{\infty}$ means the  operator norm.

\vskip 10pt\noindent {\bf Theorem 3.1.} {\it For all block-matrices
whose entries are normal matrices of same size and for all symmetric
norms,
$$
\left\|
\begin{pmatrix}
A & B \\
C & D
\end{pmatrix}
\right\| \le \left\| \, |A| + |B| + |C| + |D| \right\|.
$$
}

\vskip 10pt\noindent {\bf Theorem 3.2.} {\it For all block-matrices
whose entries are normal matrices of same size,
$$
\left\|
\begin{pmatrix}
A & B \\
C & D
\end{pmatrix}
\right\|_{\infty}
 \le \max\{\,\left\|  \, |A|+|B|\, \right\|_{\infty} ;
 \left\|  \, |C|+|D|\, \right\|_{\infty} ;
 \left\|  \, |A|+|C|\, \right\|_{\infty} ;
 \left\|  \, |B|+|D|\, \right\|_{\infty}\,\}
$$
}

\vskip 10pt\noindent {\bf Proof.}  Let $A_1,\, A_2,\ B_1,\, B_2$ be
positive and let $C_1, C_2$ be contractions. Note that
$$
A_1C_1B_1 + A_2C_2B_2=
\begin{pmatrix}
A_1 & A_2
\end{pmatrix}
\begin{pmatrix}
C_1 & 0 \\
0 & C_2
\end{pmatrix}
\begin{pmatrix}
B_1 \\ B_2.
\end{pmatrix}
$$
Applying the Cauchy-Shwarz inequality $\Vert XY\Vert \le \Vert
X^*X\Vert^{1/2} \Vert YY^*\Vert^{1/2}$ and using $\Vert ST\Vert \le
\Vert S\Vert_{\infty} \Vert T  \Vert$ show
\begin{equation*}
\left\| A_1C_1B_1 + A_2C_2B_2 \right\|\le \left\| A_1^2 +A_2^2
\right\|^{1/2} \left\| B_1^2 +B_2^2 \right\|^{1/2}.
\end{equation*}
Considering polar decompositions $A=|A^*|^{1/2}U|A|^{1/2}$ and
$B=|B^*|^{1/2}V|B|^{1/2}$ then shows that
\begin{equation}
\|A+B\|\le  \|\, |A| + |B|\, \|^{1/2}\,\|\, |A^*| + |B^*|\, \|^{1/2}
\end{equation}
for all $A,\, B$. Replacing $A$ and $B$ in (4) by
$$
\begin{pmatrix}
A & 0 \\
0 & D
\end{pmatrix}
\quad {\rm and} \quad
\begin{pmatrix}
0 & B \\ C & 0
\end{pmatrix}
$$
and using normality of $A,\, B,\, C,\, D$ yield
$$
\left\|
\begin{pmatrix}
A & B \\
C & D
\end{pmatrix}
\right\|
 \le
  \left\|
\begin{pmatrix}
|A|+|C| & 0 \\
0 & |B|+|D|
\end{pmatrix}
\right\|^{1/2}
 \left\|
\begin{pmatrix}
|A|+|B| & 0 \\
0 & |C|+|D|
\end{pmatrix}
\right\|^{1/2}. $$
 This proves Theorem 3.2. This also proves Theorem
3.1 by using the fact that
$$
\left\|
\begin{pmatrix}
X & 0 \\
0 & Y
\end{pmatrix}\right\| \le
  \left\|
X+Y \right\|
$$
for all $X,\,Y\ge0$. \qquad $\Box$

\vskip 10pt\noindent {\bf Corollary 3.3.} {\it For any partitioned
Hermitian matrix,
$$
\left\|
\begin{pmatrix}
A & X^* \\
X & B
\end{pmatrix}
\right\|_{\infty}
 \le \max\{\,\left\|  \, |A|+|X|\, \right\|_{\infty} ;
  \left\|  \, |B|+|X^*|\, \right\|_{\infty}\,\}.
$$
}

\noindent {\bf Proof.} We consider the block matrix
$$
\begin{pmatrix}
0&0&0&X \\
0&A & X^* &0\\
0&X & B   &0\\
X^*&0&0&0
\end{pmatrix}
$$
with normal blocks
$$
\begin{pmatrix}
0&0\\
0&A
\end{pmatrix} \quad
\begin{pmatrix}
B&0\\
0&0
\end{pmatrix} \quad
\begin{pmatrix}
0&X\\
X^*&0
\end{pmatrix}
$$
and we apply Theorem 3.2. \qquad $\Box$

\vskip 10pt A special case of (4) is the following statement:

\vskip 10pt\noindent {\bf Proposition 3.4.} {\it Let $A$, $B$ be
normal. Then, for all symmetric norms,
$$
\Vert A+B \Vert \le \Vert\,|A|+|B|\,\Vert.
$$
}

\noindent This can be regarded as a triangle inequality for normal
operators. In case of Hermitian operators, a stronger triangle
inequality holds [6]:

\vskip 10pt\noindent {\bf Proposition 3.5.} {\it  Let $S$, $T$ be
Hermitian. Then, for some unitaries $U,\,V$,
$$
|S+T| \le \frac{1}{2}\left\{ U(|S|+|T|)U^*+V(|S|+|T|)V^*\right\}
$$
}

\noindent  Proposition 3.5 implies Proposition 3.4 by letting
$$
S=\begin{pmatrix}
0 & A^* \\
A & 0
\end{pmatrix}
\qquad {\rm and} \qquad
T=\begin{pmatrix}
0 & B^* \\
B & 0
\end{pmatrix}.
$$
 Proposition 3.4 shows that, given $A,B\ge 0$ and any complex number $z$,
$$
\Vert A+zB \Vert \le \Vert A+|z|B\Vert
$$
so that for all integers $m=1,\, 2...,$
\begin{equation}\Vert (A+zB)^m \Vert \le \Vert (A+|z|B)^m\Vert.\end{equation}
This was observed by Bhatia and Kittaneh. They also noted the
identity
\begin{equation}
A^m+B^m = \frac{1}{m}\sum_{j=0}^{m-1} (A+wB)^j
\end{equation}
where $w$ is the primitive $m$-th root of the unit.
 Combining (5) and (6), Bhatia and Kittaneh obtained [4]: {\it Given
$A,\,B\ge0$
$$
\Vert A^m+B^m  \Vert \le \Vert (A+B)^m\Vert
$$
for all $m=1,\, 2...$ and all symmetric norms}. This result was the
starting point of superadditive or subadditive inequalities for
symmetric norms.

 \vskip 10pt {\bf References}

\noindent
{\small

\vskip 5pt\noindent [1] T.\ Ando and X.\ Zhan, Norm inequalities
related to operator monotone functions, Math.\ Ann.\ 315 (1999)
771-780. \vskip 5pt\noindent [2] J.\ S.\ Aujla and J.-C.\ Bourin,
Eigenvalue inequalities for convex and log-convex functions, Linear
Alg.\ Appl., in press (2007). \vskip 5pt\noindent [3] J.\ S.\ Aujla
and F.\ C.\ Silva, Weak majorization inequalities and convex
functions, Linear Alg.\ Appl., 369 (2003) \vskip 5pt\noindent [4]
R.\ Bhatia, F.\ Kittaneh, Norm inequalities for positive operators,
Lett.\ Math.\ Phys.\ 43 (1998) 225-231. \vskip 5pt\noindent [5]
J.-C.\ Bourin, A concavity inequality for symmetric norms, Linear
Alg.\ Appl., 413 (2006) 212-217. \vskip 5pt\noindent [6] J.-C.\
Bourin, Hermitian operators and convex functions, J.\ of Inequal.\
in Pure and Appl.\ Math.\ 6 (2005) no5, Article 139. \vskip
5pt\noindent [7] T.\ Kosem, Inequalities between $\Vert f(A+B)\Vert$
and $\Vert f(A)+f(B)\Vert$, Linear Alg.\ Appl., 418 (2006) 153-160.
\vskip 5pt\noindent [8] M.\ Uchiyama, Subadditivity of eigenvalue
sums, Proc.\ Amer.\ Math.\ Soc., 134 (2006) 1405-1412. \vskip
5pt\noindent [9] X.\ Zhan,  Matrix Inequalities, LNM 1790 (2002),
Springer, Berlin. }

\vskip 10pt
\centerline{Jean-Christophe Bourin}
 \centerline{E-mail: bourinjc@@club-internet.fr}
  \centerline{Department of mathematics,
  Kyungpook National University, Daegu 702-701, Korea}

  \vskip 10pt
\centerline{Mitsuru Uchiyama}
 \centerline{E-mail: uchiyama@@riko.shimane-u.ac.jp}
 \centerline{Department of mathematics, Interdisciplinary Faculty of Science and Engineering,}
\centerline{  Shimane University, Matsue city, Shimane 690-8504,
Japan}

\end{document}